\documentclass[a4paper,10pt]{amsart}
\usepackage{amssymb,amsmath,amsthm,hyperref}
\usepackage{diagrams}

\newcommand{\abar}{{\ensuremath{\bar{a}}}}
\newcommand{\bbar}{{\ensuremath{\bar{b}}}}
\newcommand{\cbar}{{\ensuremath{\bar{c}}}}
\newcommand{\dbar}{{\ensuremath{\bar{d}}}}

\newcommand{\xbar}{{\ensuremath{\bar{x}}}}
\newcommand{\ybar}{{\ensuremath{\bar{y}}}}
\newcommand{\zbar}{{\ensuremath{\bar{z}}}}
\newcommand{\mbar}{{\ensuremath{\bar{m}}}}
\newcommand{\Xbar}{{\ensuremath{\bar{X}}}}

\DeclareMathOperator{\tp}{tp}  

\DeclareMathOperator{\rk}{rk}

\DeclareMathOperator{\td}{td}  

\DeclareMathOperator{\loc}{Loc}   

\DeclareMathOperator{\ldim}{ldim}  
\DeclareMathOperator{\Mat}{Mat}  

\DeclareMathOperator{\ecl}{ecl} 

\newcommand{\N}{\ensuremath{\mathbb{N}}}
\newcommand{\Z}{\ensuremath{\mathbb{Z}}}
\newcommand{\Q}{\ensuremath{\mathbb{Q}}}

\newcommand{\Rexp}{\ensuremath{\mathbb{R}_{\mathrm{exp}}}}
\newcommand{\Cexp}{\ensuremath{\mathbb{C}_{\mathrm{exp}}}}

\newcommand{\Loo}{\ensuremath{L_{\omega_1,\omega}}} 
\newcommand{\Looq}{\ensuremath{\Loo(Q)}}

\newcommand{\ga}{\ensuremath{\mathbb{G}_\mathrm{a}}}   
\newcommand{\gm}{\ensuremath{\mathbb{G}_\mathrm{m}}}  

\renewcommand{\phi}{\varphi}

\renewcommand{\ge}{\ensuremath{\geqslant}}
\newcommand{\tuple}[1]{\ensuremath{\langle #1 \rangle}}
\newcommand{\class}[2]{\ensuremath{\left\{ #1 \,\left|\, #2 \right.\right\}}}

\newcommand{\iso}{\cong}

\newcommand{\into}{\hookrightarrow}

\newcommand{\subs}{\subseteq} 

\newcommand{\elsubs}{\preccurlyeq} 

\newcommand{\minus}{\ensuremath{\smallsetminus}}

\newcommand{\strong}{\ensuremath{\lhd}} 
\newcommand{\nstrong}{\ensuremath{\not\kern-4pt\lhd\;}} 

\newcommand{\gen}[1]{\ensuremath{\left\langle #1 \right\rangle}} 
\newcommand{\hull}[1]{\ensuremath{\lceil #1\rceil}}

\newcommand{\cross}{\ensuremath{\times}}

\newbox\noforkbox \newdimen\forklinewidth
\forklinewidth=0.3pt \setbox0\hbox{$\textstyle\smile$}
\setbox1\hbox to \wd0{\hfil\vrule width \forklinewidth depth-2pt
  height 10pt \hfil}
\wd1=0 cm \setbox\noforkbox\hbox{\lower 2pt\box1\lower
2pt\box0\relax}
\def\unionstick{\mathop{\copy\noforkbox}\limits}

\def\nonfork_#1{\unionstick_{\textstyle #1}}

\setbox0\hbox{$\textstyle\smile$} \setbox1\hbox to \wd0{\hfil{\sl
/\/}\hfil} \setbox2\hbox to \wd0{\hfil\vrule height 10pt depth
-2pt width
                \forklinewidth\hfil}
\newbox\doesforkbox
\setbox\doesforkbox\hbox{\lower 2pt\box1 \lower
2pt\box2\lower2pt\box0\relax}
\def\nunionstick{\mathop{\copy\doesforkbox}\limits}

\def\fork_#1{\nunionstick_{\textstyle #1}}

\newcommand{\ra}[3]{\ensuremath{#1 \stackrel{#2}{\longrightarrow} #3}}

\newcommand{\leteq}{\mathrel{\mathop:}=}

\newarrowmiddle{normal}\lhtriangle\rttriangle\uhtriangle\dhtriangle

\newarrowmiddle{iso}{\cong}{\cong}{}{}   

\newarrow{Partial}{-}{-}{-}{-}{harpoon}     

\newarrow{Equal}=====     

\newarrow{Strong}{hooka}-{normal}->     

\newarrow{Str}{}{}{normal}{}{}     

\newarrow{Iso}--{iso}->     

\newarrow{NT}====>     

\newcommand{\EAC}[1]{\ensuremath{{#1}^\sim}}  
\newcommand{\B}{\ensuremath{\mathbb{B}}}

\newcommand{\seac}{strong exponential-algebraic closedness}
\newcommand{\Seac}{Strong exponential-algebraic closedness}

\newtheorem*{prop}{Proposition}
\newtheorem{theorem}{Theorem}

\theoremstyle{definition}

\title[The axioms for Zilber's exponential fields]{A note on the axioms for Zilber's pseudo-exponential fields}
\author{Jonathan Kirby}
\date{version 2.3, \today}
\address{Jonathan Kirby\\ School of Mathematics\\ University of East Anglia\\ Norwich Research Park\\  Norwich NR4 7TJ\\ UK}
\email{jonathan.kirby@uea.ac.uk}

\begin{document}

\maketitle

\begin{abstract}
We show that Zilber's conjecture that complex exponentiation is isomorphic to his pseudo-exponentiation follows from the a priori simpler conjecture that they are elementarily equivalent. An analysis of the first-order types in pseudo-exponentiation leads to a description of the elementary embeddings, and the result that pseudo-exponential fields are precisely the models of their common first-order theory which are atomic over exponential transcendence bases. We also show that the class of all pseudo-exponential fields is an example of a non-finitary abstract elementary class, answering a question of Kes\"al\"a and Baldwin.
\end{abstract}


\section{Introduction}

Zilber's \emph{pseudo-exponential fields} were defined in \cite{Zilber05peACF0}, as models of certain axioms. Infinitary logic is used throughout the paper, which shows that the axioms are expressible by an $\Looq$-sentence. Anand Pillay emphasised to me the importance of understanding the first-order theory of these pseudo-exponential fields, which we denote by $T_\B$. In \cite{EFAI}, Zilber and I gave an axiomatization of $T_\B$, assuming the diophantine conjecture CIT. This note takes some steps towards understanding $T_\B$ unconditionally, in particular by describing the first-order types which are realised in pseudo-exponential fields. In section 2, I show the axioms for pseudo-exponential fields can be written in a first-order way, with two exceptions: an $\Loo$-sentence is needed to omit the type of a non-standard integer, and the quantifier $Q$ (there exist uncountably many) is essential for the countable closure axiom (CCP), although the axiom can be expressed in an $L(Q)$-scheme, rather than just as an \Looq-sentence. These non-elementary axioms are known to hold in $\Cexp$, so it follows immediately that if $\Cexp$ satisfies the first-order theory of pseudo-exponential fields then it satisfies the full $\Looq$-theory.

The main result of Zilber's paper is that there is a unique model of his axioms of cardinality $2^{\aleph_0}$, which we call $\B$, and indeed of each uncountable cardinal $\kappa$. Thus we have
\begin{theorem}\label{elemequiv implies isomorphism}
 If the complex exponential field \Cexp\ is elementarily equivalent to $\B$, then it is isomorphic to $\B$.
\end{theorem}
The isomorphism is Zilber's main conjecture about $\Cexp$, however Theorem~\ref{elemequiv implies isomorphism} does not make proving the isomorphism any easier, rather it shows that proving elementary equivalence is just as far out of reach. 


The class of pseudo-exponential fields naturally forms an abstract elementary class (AEC). Hyttinen and Kes\"al\"a \cite{HK06} introduced the notion of a \emph{finitary} AEC, and the class of pseudo-exponential fields has been used as an example of a non-finitary AEC, but no proof of this fact has previously appeared in the literature. From the proof that the $Q$ quantifier cannot be eliminated, it follows quickly that:
\begin{theorem}\label{not finitary}
  The category of pseudo-exponential fields (with CCP) together with closed embeddings is an abstract elementary class which is not finitary.
\end{theorem}

Section~3 of this note uses the algebra of exponential fields developed in \cite{FPEF} and \cite{EFAI} to address the issues of elementary embeddings and atomicity for pseudo-exponential fields. In fact here we can work with a broader class, not assuming the countable closure axiom, but unfortunately still narrower than the class of all models of $T_\B$. The main technical difficulty is that, without CIT, we do not understand what saturated models of $T_\B$ look like, and we do not have an unconditional quantifier-elimination result, so the proofs have to use less direct techniques including a careful analysis of the types which are realised.

\begin{theorem}\label{elem embeddings}
The strong embeddings between pseudo-exponential fields (not necessarily with CCP) are exactly the elementary embeddings. In particular, all pseudo-exponential fields, including those of finite exponential transcendence degree, are elementarily equivalent.
\end{theorem}
This theorem complements Theorem~5.13 of \cite{Zilber05peACF0}, which adds the additional hypothesis of infinite exponential transcendence degree but gets the stronger conclusion that strong embeddings are \Loo-embeddings. The stronger conclusion fails in general. 

Finally we show:
\begin{theorem}\label{atomic models}
  Each pseudo-exponential field (not necessarily with CCP) is an atomic model of $T_\B$ over an exponential transcendence base, in the language $\tuple{+,\cdot,\exp}$. In particular, $T_\B$ has a prime model. Conversely, every model of $T_\B$ which is atomic over an exponential transcendence base is a pseudo-exponential field.
\end{theorem}

Zilber states in \cite[Theorem~5.16]{Zilber05peACF0} that the uncountable pseudo-exponential fields (with the countable closure property) are prime over exponential transcendence bases but he uses an expanded language $L^*$ which is an expansion of our language by $\Loo$-definitions, not by first-order definitions.  So Theorem~\ref{atomic models} does not follow directly from Zilber's methods.

Some of the work for this paper was done while I was visiting the University of Helsinki, and I am grateful for their hospitality.

\section{The axioms for pseudo-exponentiation}
We give five axioms capturing Zilber's definition. Explanations of the terminology in axioms 4 and 5 are deferred to the more detailed discussions afterwards. We use only the language $\tuple{+,\cdot,\exp}$.
\begin{description}
 \item[1. ELA-field] $F$ is an algebraically closed field of characteristic zero, and its exponential map $\exp$ is a homomorphism from its additive group to its multiplicative group, which is surjective.

 \item[2. Standard kernel] the kernel of the exponential map is an infinite cyclic group generated by a transcendental element $\tau$.

 \item[3. Schanuel Property] The \emph{predimension function} 
\[\delta(\xbar) \leteq \td(\xbar, \exp(\xbar))- \ldim_\Q(\xbar)\]
satisfies $\delta(\xbar) \ge 0$ for all tuples $\xbar$ from $F$.

\item[4. \Seac] If $V$ is a rotund, additively and multiplicatively free subvariety of $\ga^n\cross \gm^n$ defined over $F$ and of dimension $n$, and $\abar$ is a finite tuple from $F$, then there is $\xbar$ in $F$ such that $(\xbar,e^\xbar) \in V$ and is generic in $V$ over $\abar$.

\item[5. Countable Closure Property] For each finite subset $C$ of $F$, the exponential algebraic closure $\ecl^F(C)$ of $C$ in $F$ is countable. 
\end{description}

Clearly axiom 1 is first-order expressible.

\subsection{Standard Kernel}

Write $\ker(F) = \class{x\in F}{\exp(x) = 1}$, the kernel of the exponential map.
Define 
\[Z(F) = \class{r \in F}{\forall x[x \in \ker \to rx \in \ker]},\]
the \emph{multiplicative stabilizer} of the kernel. 

Axiom 2 can be split into three parts.
\begin{description}
 \item[2a] The kernel is a cyclic $Z$-module.
 \item[2b] Every element of the kernel is transcendental over $Z$.
 \item[2c Standard integers] 
 \[(\forall r \in Z) \bigvee_{n \in \N} r = \underbrace{1+\cdots+1}_n \quad \vee\quad r + \underbrace{1+\cdots+1}_n = 0\]
\end{description}
 It is clear that axioms 2a and 2b are first-order expressible, and that 2c is not, but it is given by a single $\Loo$-sentence. Equivalently, 2c can be viewed as omitting the incomplete type of a nonstandard integer.

\subsection{The Schanuel Property}
Axiom 3, the Schanuel Property, is equivalent to the following axiom scheme, which is first-order expressible provided axiom 2 holds so we can quantify over the standard integers:
for each $n \in \N$, for each subvariety $V \subs \ga^n \cross \gm^n$ defined over $\Q$, of dimension $n-1$,
\[(\forall x_1,\ldots,x_n)(\exists \bar{m} \in \Z^n \minus \{\bar{0}\})\left[(\bar{x},\exp(\bar{x})) \in V \to \sum_{i=1}^n m_i x_i = 0\right].\]

\subsection{\Seac}

Write $G$ for the algebraic group $\ga \cross \gm$. Each matrix $M \in \Mat_{n \cross n}(\Z)$ defines a homomorphism
$\ra{G^n}{M}{G^n}$ by acting as a linear map on $\ga^n$ and as a
multiplicative map on $\gm^n$. If $V \subs G^n$, we write $M\cdot V$
for its image. Note that if $V$ is a subvariety of $G^n$, then so is $M\cdot
V$.

  An irreducible subvariety $V$ of $G^n$ is said to be \emph{rotund} iff for
  every matrix $M \in \Mat_{n \cross n}(\Z)$ we have $\dim M\cdot V \ge \rk M$.

Suppose that $(\xbar,\ybar)$ is a generic point of $V$ over $F$, with the $x_i$ being the coordinates from $\ga$ and the $y_i$ the coordinates from $\gm$. We say that $V$ is \emph{multiplicatively free} iff the $y_i$ do not satisfy any equation of the form  $\prod_{i=1}^n y_i^{m_i} = b$ with the $m_i \in \Z$, not all zero, and $b \in \gm(F)$. Equivalently, the projection of $V$ to $\gm^n$ does not lie in any coset of a proper algebraic subgroup of $\gm^n$. Similarly, we say that $V$ is \emph{additively free} iff the $x_i$ do not satisfy any equation of the form $\sum_{i=1}^n m_i x_i = a$ with the $m_i \in \Z$, not all zero, and $a \in F$.

\begin{prop}\label{SEAC fo}
  Axiom 4, \seac, is first-order expressible modulo axioms 1, 2, and 3.
\end{prop}
\begin{proof}
We consider parametric families $(V_p)_{p \in P}$ of subvarieties of $G^n$, where $P$ is some parametrizing variety. It is a well-known fact (part of the \emph{fibre dimension theorem}) that the set of $p$ such that $V_p$ is irreducible and of dimension $n$ is first-order definable in the field language. The property of being additively free is not definable in the field language, since for example the subvariety of $G^2$ given by the equation
\[x_1 + p x_2 = 0\]
is additively free iff $p \notin \Q$. However, it is definable as follows allowing quantification over $\Z$:
\[(\forall \mbar \in \Z^n\minus \{\bar{0}\})\forall z\exists \xbar\left[\xbar \in V_p \wedge \sum_{i=1}^n m_i x_i \neq z\right].\]
It is easy to give similar definitions showing that rotundity and multiplicative freeness are definable allowing quantification over $\Z$. However, Theorem~3.2 of \cite{Zilber05peACF0} shows that these two properties are even first-order definable in the field language. For a parametric family $(V_p)_{p \in P}$ of subvarieties of $G^n$, let $P'$ be the set of $p \in P$ consisting of those $p$ such that $V_p$ is irreducible, of dimension $n$, rotund, and additively and multiplicatively free. Consider the following axiom scheme with one axiom for each family $(V_p)_{p \in P}$ and each natural number $r$:
\begin{multline*}
  (\forall p\in P')(\forall \abar \in F^r)(\exists \xbar \in F^n)(\forall \mbar \in \Q^{n+r}) \\
\left[(\xbar,e^{\xbar}) \in V_p \wedge \left(\sum_{i=1}^n m_i x_i + \sum_{i=1}^r m_{n+i} a_i = 0 \to \bigwedge_{i=1}^n m_i = 0\right)\right]
\end{multline*}
This scheme is first-order expressible assuming axiom 2 holds, so we can quantify over $\Z$ and hence over $\Q$. It also follows from axiom 4, since that axiom gives an $\xbar$ such that $(\xbar,e^\xbar)$ is generic in $V_p$ over $\abar$. Since $V_p$ is additively free, that $\xbar$ does not satisfy any $\Q$-linear equation over $\abar$, which is all the scheme requires.

Now let $V$ be any rotund, additively and multiplicatively free subvariety of $G^n$ defined over $F$ and of dimension $n$, and $\abar$ be a finite tuple from $F$, as in the hypothesis of axiom 4. Then $V$ is $V_p$ for some parametric family $(V_p)_{p \in P}$ as above, with $p \in P'$. Using axiom 3, by extending $\abar$ if necessary, we may assume that $\delta(\ybar/\abar) \ge 0$ for all tuples $\ybar$ from $F$, and also that $V$ is defined over $\abar$. Now invoke the axiom scheme to find an $\xbar$. We have $\ldim_\Q(\xbar/\abar) = n$, so we must have $\td(\xbar,e^{\xbar}/\abar,e^{\abar}) \ge n$. But $\dim V = n$, so $(\xbar,e^\xbar)$ must be generic in $V$ over $(\abar,e^\abar)$, so a fortiori over $\abar$. Thus axiom 4 is equivalent to this scheme, modulo axioms 1, 2, and 3.
 \end{proof}

On page 87 of \cite{Zilber05peACF0}, Zilber remarks:
\begin{quotation}
\noindent The definition [of strongly exponentially-algebraically closed] assumes a ``slight saturatedness'' of the exponentially-algebraically closed structure.
\end{quotation}
This remark had led me to assume that \seac\ was not first-order, even assuming the other axioms, so the above result was somewhat unexpected. Indeed the fact that \seac\ is first-order means that the Zilber's other notion of \emph{exponential-algebraic closedness} is redundant for the construction of the exponential fields which do have standard kernel and the Schanuel property (although it is used in this note, in \S\ref{et type} where both these properties fail). In  \cite[Theorem~5.5]{EFAI} we prove that the notions of exponential-algebraic closedness and \seac\ coincide under the additional assumption that the diophantine conjecture CIT is true (and assuming the other relevant axioms) and it would be interesting to know if the same result can be proved unconditionally.

\subsection{The countable closure property}\label{ecl section}

In any exponential field $F$ there is a pregeometry called exponential algebraic closure, which we write $\ecl^F$. We give a quick account of its definition. Details can be found in  \cite{Macintyre96} or \cite{EAEF}.  
An exponential polynomial (without iterations of exponentiation) is a function of the form $f(\Xbar) = p(\Xbar,e^\Xbar)$ where $p \in F[X_1,\ldots,X_n,Y_1,\ldots,Y_n]$ is a polynomial. We can extend the formal differentiation of polynomials to exponential polynomials in a unique way such that $\frac{\partial e^X}{\partial X} = e^X$.

  A \emph{Khovanskii system of width $n$} consists of exponential
  polynomials $f_1,\ldots,f_n$ with equations
  \begin{equation}\label{Khov1}
f_i(x_1,\ldots,x_n) = 0 \quad \mbox{for } i=1,\ldots,n
\end{equation}
and the inequation
\begin{equation}\label{Khov2}
\begin{vmatrix}
  \frac{\partial f_1}{\partial X_1} & \cdots &\frac{\partial
    f_1}{\partial X_n}\\
  \vdots & \ddots & \vdots \\
  \frac{\partial f_n}{\partial X_1} & \cdots &\frac{\partial
    f_n}{\partial X_n} \end{vmatrix} (x_1,\ldots,x_n) \neq 0.
    \end{equation}
where the differentiation here is the formal differentiation of exponential polynomials.

  For any subset $C$ of $F$, we define $a \in \ecl^F(C)$ iff
  there are $n \in \N$, $a_1,\ldots,a_n \in F$, and exponential polynomials $f_1,\ldots,f_n$
  with coefficients from $\Q(C)$ such that $a = a_1$ and
  $(a_1,\ldots,a_n)$ is a solution to the Khovanskii system given by
  the $f_i$.

  We say that $\ecl^F(C)$ is the \emph{exponential algebraic closure} of $C$ in $F$. If $a \in \ecl^F(C)$ we say that $a$ is \emph{exponentially algebraic} over $C$ in $F$, and otherwise that it is \emph{exponentially transcendental} over $C$ in $F$.

Theorem~1.1 of \cite{EAEF} states that $\ecl^F$ is a pregeometry in any exponential field. We have stated axiom 5 in terms of this pregeometry. The definition of the pregeometry originally used by Zilber is different, and makes sense only assuming axiom 2. However, when axiom 2 holds, the two definitions agree by \cite[Theorem~1.3]{EAEF}. With Zilber's original definition one can see that axiom 5 is expressible as an \Looq-sentence. Using this definition we can show:
\begin{prop}
Axiom 5, the countable closure property, is expressible as an $L(Q)$-scheme.
\end{prop}
\begin{proof}
For the Khovanskii system on exponential polynomials $f = (f_1,\ldots,f_n)$, write $\chi_f(\xbar,\zbar)$ for the first-order formula expressing (\ref{Khov1}) and (\ref{Khov2}), where $\zbar$ denotes the coefficients for the exponential polynomials.
Then axiom 5 is expressed by the L(Q)-sentences
\[(\forall \zbar)\neg(Q x_1)(\exists x_2,\ldots,x_n)\chi_f(\xbar,\zbar)\]
where $f$ ranges over all finite lists of exponential polynomials with variables $\zbar$ as coefficients.
\end{proof}

\subsection{The complex exponential field and proof of Theorem~\ref{elemequiv implies isomorphism}}

The complex exponential field $\Cexp$ is the field of complex numbers equipped with the usual complex exponential function given by $\exp(z) = \sum_{n \in \N} \frac{z^n}{n!}$. Axioms 1 and 2 are chosen such that \Cexp\ satisfies them. Zilber noted \cite[Lemma~5.12]{Zilber05peACF0} that \Cexp\ also satisfies axiom 5, the countable closure property. With the definition of $\ecl^F$ from \S\ref{ecl section}, we can give a shorter proof. Given a finite subset $C$ of \Cexp, there are only countably many Khovanskii systems with coefficients from $\Q(C)$. The inequation in a Khovanskii system says that the Jacobian of the functions $f_1,\ldots,f_n$ does not vanish so, by the implicit function theorem, solutions to a Khovanskii system are isolated in the complex topology. Hence there are only countably many solutions to each system, so $\ecl^{\Cexp}(C)$ is countable.

We have seen that axioms 3 and 4 are first-order expressible modulo axioms 1 and 2, and certainly \Cexp\ has cardinality $2^{\aleph_0}$, so Theorem~\ref{elemequiv implies isomorphism} follows.

\subsection{Strong extensions}\label{strong extensions}

We now summarize the definitions and results from \cite{FPEF} and \cite{EFAI} which we shall need.

In this note, a \emph{partial exponential subfield} of an ELA-field $F$ is a subfield $F_0 \subs F$ together with a $\Q$-linear subspace $D(F_0)$ of $F_0$ and the restriction of the exponential map $\exp_F$ to $D(F_0)$, such that $F_0$ is generated as a subfield of $F$ by $D(F_0) \cup \exp(D(F_0))$. Thus $F_0$ is determined by $D(F_0)$. For any finite tuple $\abar$ from $F$ we define
\[ \delta(\abar/F_0) = \td(\abar,e^\abar/F_0) - \ldim_\Q(\abar/D(F_0))\]
and say $F_0$ is strong in $F$, written $F_0 \strong F$, iff for all $\abar \in F$, $\delta(\abar/F_0) \ge 0$.

If $B \subs F$ is a subset, then $\gen{B}_F$, the partial exponential subfield of $F$ generated by $B$, is the partial exponential subfield $F_0$ with $D(F_0)$ equal to the $\Q$-linear span of $B$. We write $B \strong F$ iff $\gen{B}_F \strong F$.

If $F$ satisfies the Schanuel property or, more generally, if there is $C \subs B$ with $C \strong F$, then there is a smallest partial exponential subfield of $F$ which contains $B$ and is strong in $F$. We write it as $\hull{B}_F$ and call it the \emph{hull} of $B$ in $F$. For any subset $B \subs F$, we write $\gen{B}^{ELA}_F$ for the smallest ELA-subfield of $F$ containing $\ker(F) \cup B$, and $\hull{B}^{ELA}_F$ for $\gen{\hull{B}_F}^{ELA}_F$.

\begin{prop}[A]
Suppose $F$ is an ELA-field with standard kernel, and $\bbar$ is a finite tuple from $F$ such that $\bbar \strong F$. Then there is $m \in \N^+$ such that the isomorphism type of $\gen{\bbar}^{ELA}_F$ is determined by the algebraic locus $V = \loc(\bbar/m,e^{\bbar/m}/\ker(F))$, a subvariety of $G^n$, and does not depend on $F$. Furthermore, $\gen{\bbar}^{ELA}_F$ is strong in $F$. 
\end{prop}
\begin{proof}
The existence of some $m \in \N^+$ such that the partial exponential subfield $\gen{\bbar}_F$ is determined up to isomorphism by $V$ follows from the Thumback Lemma, \cite[Fact~2.15]{FPEF}, or see \cite[Theorem~2]{Zilber06covers} for a proof. Then \cite[Theorem~2.18]{FPEF} applies to show that $\gen{\bbar}^{ELA}_F$ is uniquely determined and strong in $F$.
\end{proof}

\begin{prop}[B]
If $F \strong M$ is a strong extension of ELA-fields such that $\ker(F) = \ker(M)$, $M$ is generated as an ELA-field by $F$ and a finite tuple $\bbar \in M^n$, and either $F$ is countable or $\ker(F)$ is $\aleph_0$-saturated, then there is $m \in \N^+$ such that the isomorphism type of $M$ as an extension of $F$ is determined by $V = \loc(\bbar/m,e^{\bbar/m}/F)$. Furthermore, $\bbar$ can be chosen to be $\Q$-linearly independent over $F$, and then $V$ is rotund, additively and multiplicatively free, and of dimension at least $n$. The extension is exponentially algebraic, that is $M = \ecl^{M}(F)$, iff $\dim(V) = n$.
\end{prop}
In this case, we write $M$ as $F|V$, ``$F$ extended by $V$''.
\begin{proof}
See \S3 and \S5 of \cite{FPEF} for the case where $F$ is countable. In the case where $\ker(F)$ is $\aleph_0$-saturated, the analysis is exactly the same except that one needs to use \cite[Theorem~3.3]{EFAI} in place of \cite[Theorem~2.18]{FPEF}.
\end{proof}

\subsection{\Loo-theory}
 Let $\Psi$ be an \Loo-sentence expressing axioms 1---4. For any natural number $n$, it is easy to give an \Loo-sentence $\Phi_n$ specifying that $F$ has exponential-transcendence degree equal to $n$. 
 
From \S\ref{strong extensions}, we see that axiom 4, \seac, is equivalent (assuming axioms 1, 2, and 3) to an existential closedness property for exponentially-algebraic strong ELA-extensions, which do not extend the kernel. Lemma~5.9 of \cite{FPEF} shows these extensions can be freely amalgamated, so, at least in the countable case, we can characterize the models as certain Fra\"iss\'e limits. Indeed any ELA-field $F$ which is countable or has $\aleph_0$-saturated kernel has a unique smallest strongly exponentially-algebraically closed extension, which we write as $\EAC{F}$. For more details, see \cite[\S6]{FPEF}.

In particular, each of the \Loo-sentences $\Psi_n \leteq \Psi \wedge \Phi_n$ and $\Psi_\infty \leteq \Psi \wedge \bigwedge_{n \in \N} \neg \Phi_n$ is countably categorical by the uniqueness of Fra\"iss\'e limits (or specifically by \cite[Corollary~6.10]{FPEF}), and hence complete. For $\Psi_\infty$ this was already proved in \cite[Theorem~5.13]{Zilber05peACF0}. Clearly these are the only completions of $\Psi$ as \Loo-theories. In particular, $\Psi_\infty$ gives the complete \Loo-theory of $\B$. We write $B_n$ for the unique countable model of $\Psi_n$.

\subsection{Necessity of Q and proof of Theorem~\ref{not finitary}}

We now show that the $Q$ quantifier cannot be eliminated.
\begin{prop}
  The countable closure property is not expressible in $L_{\infty,\omega}$, even modulo axioms 1---4.
\end{prop}
\begin{proof}
 The idea is to find and use an exponentially algebraic regular (in fact strongly minimal) type which is totally categorical and orthogonal to the generic (exponentially transcendental) type, and to the kernel.

 Let $F_0 = B_0$, and adjoin an extra element $a$ such that $e^a = a$. This $a$ generates a well-defined ELA-field extension $F_0|V$ of $F_0$, where $V$ is the subvariety of $G$ given by the equation $x = y$. Let $F_1 = (F_0|V)^\sim$, the strong exponential-algebraic closure of $F_0|V$. By construction, $F_1$ satisfies axioms 1---5, and is an exponentially-algebraic extension of $F_0$, so has the same exponential transcendence degree, $0$. Thus $F_1 \iso F_0$. 

Now suppose that $\bbar$ is any finite tuple from $F_0$ and let $A = \hull{\bbar}^{ELA}_{F_0}$. Then $A \strong F_0$ and since $F_0 \strong F_1$ we also have $A \strong F_1$. By the general Fra\"iss\'e construction, or more specifically by Proposition~6.9 of \cite{FPEF}, $F_0$ and $F_1$ are both isomorphic over $A$ to its strong exponential-algebraic closure $A^\sim$, so in particular to each other. The same holds for any finite tuple \bbar, so it follows that the inclusion $F_0 \into F_1$ is an $L_{\infty,\omega}$-embedding. Now we can iterate the construction to get a chain
\[F_0 \into F_1 \into F_2 \into \cdots \into F_\alpha \into \cdots\]
of length $\omega_1$. The union $F_{\omega_1}$ of this chain will be $L_{\infty,\omega}$-equivalent to each element of the chain, hence to $B_0$, but $\class{a \in F_{\omega_1}}{e^a=a}$ has cardinality $\aleph_1$, and each such element is exponentially algebraic, so the countable closure property fails.
\end{proof}

We say that an embedding $F_1 \into F_2$ of pseudo-exponential fields is a \emph{closed embedding} iff the image of $F_1$ is exponentially-algebraically closed in $F_2$, that is, $\ecl^{F_2}(F_1) = F_1$. The category of all pseudo-exponential fields together with closed embeddings forms an abstract elementary class (AEC). The notion of an AEC being \emph{finitary} was introduced by Hyttinen and Kes\"al\"a \cite{HK06} and studied also by Kueker \cite{Kueker08}. Using Definition~3.1 of Kueker's paper, an AEC $\mathcal{K}$ is finitary iff whenever $M, N \in \mathcal{K}$ and $f: M \into N$ is an L-embedding such that, for every finite tuple $\abar$ from $M$, there is a $\mathcal{K}$-embedding $g : M \to N$ such that $g(\abar) = f(\abar)$, then $f$ is a $\mathcal{K}$-embedding.

In the above proof we showed that the inclusion map $F_0 \into F_1$ satisfied this last condition, but clearly it is not a closed embedding. Thus Theorem~\ref{not finitary} is established. The Proposition shows that our AEC is not closed under  \Loo-equivalence, and one could use Kueker's result that finitary AECs are closed under  \Loo-equivalence to give another, less direct, proof of Theorem~\ref{not finitary}.

\section{The first-order theory and elementary embeddings}

\subsection{The exponentially transcendental type}\label{et type}

Let $f_1,\ldots,f_n$ be exponential polynomials in variables $x_1,\ldots,x_n$, and, as before, write $\chi_f(\xbar)$ for the first-order formula corresponding to the Khovanskii system on the $f_i$, now suppressing the variables $\zbar$ corresponding to coefficients of the $f_i$.

Given any exponential field $F$, and any set $A$ of parameters from $F$, the exponentially transcendental type over $A$ is the set of formulas 
\[ (\forall x_2,\ldots,x_n)[\neg \chi_f(x,x_2,\ldots,x_n)]\]
where $n$ ranges over all positive natural numbers and $f$ ranges over all $n$-tuples of exponential polynomials with coefficients from $A$. Write $p|_A(x)$ for this type over $A$. For the exponential fields under consideration it is a consistent partial type. In some cases, for example $\Rexp$, the type is not complete.

\begin{prop}
In $T_\B$, for any set of parameters $A$, the type $p|_A(x)$ is complete.
\end{prop}
\begin{proof}
It is enough to prove the result for finite $A$. So let $M$ be an $\aleph_0$-saturated model of $T_\B$, and $A$ a finite subset of $M$. Suppose $a,b \in M$ are each exponentially transcendental over $A$. Let $M_0 = \ecl^M(A)$. We will set up a back-and-forth system showing that $\tp(a/M_0) = \tp(b/M_0)$. In particular $\tp(a/A) = \tp(b/A)$.

Firstly, note that the partial E-field extensions of $M_0$ generated by $a$ and $b$ are isomorphic, since we have $\td(a,e^a/M_0) = \td(b,e^b/M_0) = 2$. Furthermore, since $a$ and $b$ are exponentially transcendental over $M_0$, these partial E-fields are strong in $M$. By Proposition~\ref{strong extensions}(B), there is an isomorphism 
\[\theta_1: \gen{M_0,a}^{ELA}_M \rIso \gen{M_0,b}^{ELA}_M\]
between the ELA-subfields of $M$ generated by $M_0 \cup\{a\}$, and $M_0 \cup \{b\}$, fixing $M_0$ and sending $a$ to $b$. Furthermore, these ELA-subfields of $M$ are strong in $M$.

Now suppose we have $n$-tuples $\abar$ and $\bbar$ in $M$, each $\Q$-linearly independent over $M_0$, with $a_1 = a$ and $b_1=b$ and an isomorphism
\[\theta_n: \hull{M_0,\abar}^{ELA}_M \rIso \hull{M_0,\bbar}^{ELA}_M\]
between the strong ELA-subfields of $M$ generated by  $M_0 \cup \abar$, and $M_0 \cup \bbar$,
fixing $M_0$ and sending $a_i$ to $b_i$ for $i=1,\ldots,n$.

Let $c \in M$. We want to find $d \in M$ and an isomorphism
\[\theta_{n+1}: \hull{M_0,\abar,c}^{ELA}_M \rIso \hull{M_0,\bbar,d}^{ELA}_M\]
extending $\theta_n$. Write $F = \hull{M_0,\abar}^{ELA}_M$. By extending $\abar$ and $\bbar$, we may assume $M_0 \cup \abar \strong M$, so $F = \gen{M_0,\abar}^{ELA}_M$, and similarly $F' \leteq \gen{M_0,\bbar}^{ELA}_M \strong M$.

If $c \in F$ then take $d= \theta_n(c)$. If $c$ is exponentially transcendental over $F$ then, using $\aleph_0$-saturation, we can find $d$ exponentially transcendental over $A \cup \bbar$, which is therefore exponentially transcendental over $F' = \hull{M_0,\bbar}^{ELA}_M$, and the same argument as above gives us $\theta_{n+1}$.

Otherwise, since $F \strong M$, there is a finite tuple $\cbar$ of shortest length extending $c$ such that $\delta(\cbar/F) = 0$. By Proposition~\ref{strong extensions}(B), $\cbar$ generates an ELA-extension of the form $F \strong F|V$ where $V$ is the locus of $(\cbar,e^\cbar)$ over $F$, and $V$ is  rotund, additively and multiplicatively free, and of dimension equal to the length of the tuple $\cbar$. Let $\alpha$ be a finite tuple of parameters from $F$ over which $V$ is defined, and consider the variety $V'$ defined by the same formula as $V$ but with parameters $\beta = \theta_n(\alpha)$. Without loss of generality we may assume that $\bbar$ is $\Q$-linearly independent over $M_0$. 

$M$ was chosen to be an $\aleph_0$-saturated model of $T_\B$. As such, we do not know that $M$ is strongly exponentially-algebraically closed, because we do not know that axiom 4 is unconditionally first-order. However, $M$ does satisfy the axiom of exponential-algebraic closedness, (see \cite[\S5.5]{EFAI}), which in combination with the $\aleph_0$-saturation gives us $(\dbar,e^\dbar) \in V'(M)$, generic in $V'$ over $\beta \cup \bbar$, and such that no $\Q$-linear combination of $\dbar$ and $\bbar$ is exponentially algebraic over $A$, because these conditions can be expressed by a partial type over the finite set $\beta \cup \bbar \cup A$.

Then $\dbar$ is $\Q$-linearly independent over $M_0 \cup \bbar$, so $(\dbar,e^\dbar)$ is generic in $V'$ over 
$M_0\cup\bbar\cup e^\bbar$, and, using Proposition~\ref{strong extensions}(B) once more, $(\dbar,e^\dbar)$ is generic in $V'$ over $F'$. Hence the ELA-extension of $F'$ generated by $\dbar$ is of the form $F'|V'$, and so we have an isomorphism $\theta_{n+1}$ extending $\theta_n$ as required.
\end{proof}

\subsection{Non-isolation of the exponentially transcendental type}\label{nonisolation}

\begin{prop}
If $M \models T_\B$ and $A$ is any set of parameters from $M$, then the type $p|_A(x)$ is not isolated.
\end{prop}

\begin{proof}
By completeness of $p|_A(x)$ and the compactness theorem, if it were isolated by some formula then it would be isolated by a finite subtype of the form 
\[q(x,\abar) = \bigwedge_{i=1}^r(\forall x_2,\ldots,x_{n_i})[\neg\chi_{f_i}(x,x_2,\ldots,x_{n_i},\abar)]\]
where $\abar$ is the finite tuple from $A$ consisting of all coefficients from the exponential polynomials in the $f_i$. Let $N = \max\class{n_i}{i=1,\ldots,r}$.

Define $e_0(x) = x$ and $e_{n+1}(x) = \exp(e_n(x))$ for $n \in \N$. Then in fact for any tuple $\bbar$ in \B, there is $c \in \B$ such that $e_{N+1}(c) = c$ and $c,e_1(c),\ldots,e_N(c)$ are algebraically independent over $\bbar$, by \seac. Then $\B \models q(c,\bbar)$ because any tuple witnessing that $c$ is exponentially algebraic must span the $N+1$-dimensional $\Q$-vector space spanned by $\class{e_n(c)}{n = 0,\ldots,N}$, but the Khovanskii systems in $q$ look only at tuples of length up to $N$. Hence 
\[T_\B \vdash \forall \ybar \exists x [q(x,\ybar) \wedge e_{N+1}(x) = x].\]
 So there is $c \in M$ with $M \models q(c,\abar)$ but $c$ exponentially algebraic over $\abar$. Hence $p|_A(x)$ is not isolated.
\end{proof}

\subsection{Isolated types}\label{isolation}

Let $F \models \Psi$ and let $X$ be an exponential transcendence base for $F$. Any element of $\Z$ lies in the definable closure of $\emptyset$, and the kernel generator $\tau$ is model-theoretically algebraic, with $-\tau$ being the only conjugate. So the kernel of $F$ is model-theoretically algebraic. Now let $\abar$ be an $n$-tuple from $F$, and let $X_0$ be the smallest subset of $X$ such that $\abar \in \ecl(X_0)$. Let $\bbar$ be an $r$-tuple which is a $\Q$-linear basis for $D(\hull{X_0, \tau, \abar}_F)$ over $X_0 \cup \{\tau\}$, let $V = \loc(\bbar,e^\bbar/X_0,\tau)$, and let $M \in \Mat_{n \times r}(\Q)$ be such that $\abar = M \bbar$. Consider the first-order formula $\phi(\xbar)$ given by
\[\exists \ybar [(\ybar,e^\ybar) \in V \wedge \ybar \mbox{ is $Q$-linearly independent over $X_0$} \wedge \xbar = M\ybar]\]
where here $Q$ is the field of fractions of the definable subring $Z$, not the quantifier. Since $F$ has standard integers, $Q(F)$ will be $\Q$.
Now if $F \models \phi(\abar')$ and $\bbar'$ is a witness for $\ybar$ then $\ldim_\Q(\bbar'/X_0,\tau) = r$, but $\dim V = r$ and $\delta(\bbar'/X_0) = 0$, so $(\bbar,e^\bbar)$ is generic in $V$ over $X_0 \cup \{\tau\}$. Hence $\bbar' \in \ecl(X_0)$ and, by Proposition~\ref{strong extensions}(A), there is an isomorphism between $\gen{X_0,\bbar}^{ELA}_F$ and $\gen{X_0,\bbar'}^{ELA}_F$, preserving $X_0$ and $\tau$ and sending $\bbar$ to $\bbar'$. Since $X \minus X_0$ is exponentially-algebraically independent over $\ecl^F(X_0)$, this extends to an isomorphism between  $\gen{X_1,\bbar}^{ELA}_F$ and $\gen{X_1,\bbar'}^{ELA}_F$ for any countable subset $X_1$ of $X$. Then by \seac\ of $F$, $\bbar$ and $\bbar'$ are back-and-forth equivalent over any such $X_1$, and hence have the same first-order type over $X$. Hence $\abar$ and $\abar'$ have the same type over $X$. Thus the formula $\phi$ isolates the first-order type of $\abar$ over $X$.

\subsection{Elementary equivalence}
We can now prove the second part of Theorem~\ref{elem embeddings}, showing that all the $B_n$ for $n \in \N$ are elementarily equivalent to $\B$. Let $\abar$ be an exponentially algebraically independent $n$-tuple in \B, and add parameters for \abar. The type of a non-standard integer is a non-isolated (partial) type, and, by the above, the exponentially transcendental type over \abar\ is also non-isolated. Hence, by the omitting types theorem, there is a countable model $M$ of $T_\B$ containing $\abar$ which omits both types. Then $M \models \Psi_n$ and, by countable categoricity of $\Psi_n$, $M$ is isomorphic to $B_n$. Hence $B_n \models T_\B$.

\subsection{Atomic models and proof of Theorem~\ref{atomic models}}
From \ref{isolation} above, the first-order type over an exponential transcendence base of any finite tuple from any model $F \models \Psi$ is isolated by a single first-order formula in the language $\tuple{+,\cdot,\exp}$, so each such $F$ is an atomic model of the expansion of $T_\B$ by parameters for the exponential transcendence base. In particular, $B_0$ is atomic over $\emptyset$, and a countable atomic model is prime, so $B_0$ is the prime model of $T_\B$. For the converse statement, any model of $T_\B$ which is atomic over an exponential transcendence base has standard integers and hence is a model of $\Psi$, that is, of axioms 1 -- 4. That completes the proof of Theorem~\ref{atomic models}.

\subsection{Elementary embeddings and the end of the proof of  Theorem~\ref{elem embeddings}}
If an embedding $F_1 \subs F_2$ of models of $\Psi$ is not strong then there are finite tuples $\abar \strong F_1$ and $\bbar \in F_2^n$ such that $\bbar$ is $\Q$-linearly independent over $\abar$ and $\delta(\bbar/\abar) <0$, so $(\bbar,e^\bbar)$ lies in some algebraic variety $V$ defined over $(\abar,e^\abar)$ of dimension less than $n$. The first-order formula expressing the existence of such a $\bbar$ is true of $\abar$ in $F_2$ but false in $F_1$, and hence the inclusion of $F_1$ in $F_2$ is not elementary.

Now suppose $F_1 \strong F_2$ is a strong extension of models of $\Psi$. We use the Tarski-Vaught test to show that $F_1 \elsubs F_2$. So let $\abar$ be a tuple from $F_1$, $b\in F_2$, and $\phi(\xbar,y)$ a first-order formula such that $F_2 \models \phi(\abar,b)$. Extending $\abar$ if necessary, we may assume that $\abar \strong F_1$, and so $\abar \strong F_2$. Since the kernel does not extend, the ELA-subfield $F_0$ of $F_2$ generated by $\abar$ is contained in $F_1$. If $b$ is exponentially algebraic over $\abar$, then it generates a strong ELA-extension of the form $F_0 \strong F_0 | V$ for some perfectly rotund $V$. There is an isomorphic strong ELA-extension inside $F_1$ by \seac\ of $F_1$, and the element $c$ corresponding to $b$ under the isomorphism realises the same principal formula as $b$ in $F_2$. In particular, $F_2 \models \phi(\abar,c)$. Otherwise $b$ is exponentially transcendental over $\abar$, and, since the exponentially transcendental type is approximated by isolated exponentially algebraic types as in \ref{nonisolation} above, we can again find such a $c$ in $F_1$, depending on the formula $\phi$. So $F_1 \elsubs F_2$. That completes the proof of Theorem~\ref{elem embeddings}.


\begin{thebibliography}{Mac96}

\bibitem[HK06]{HK06}
T.~Hyttinen and Meeri Kes{\"a}l{\"a}.
\newblock Independence in finitary abstract elementary classes.
\newblock {\em Ann. Pure Appl. Logic}, 143(1-3):103--138, 2006.

\bibitem[Kir10]{EAEF}
Jonathan Kirby.
\newblock Exponential {A}lgebraicity in {E}xponential {F}ields.
\newblock {\em Bull. London Math. Soc.}, 42(5):879--890, 2010.

\bibitem[Kir11]{FPEF}
Jonathan Kirby.
\newblock Finitely presented exponential fields.
\newblock Submitted, arXiv:0912.4019v4, July 2011.

\bibitem[Kue08]{Kueker08}
David~W. Kueker.
\newblock Abstract elementary classes and infinitary logics.
\newblock {\em Ann. Pure Appl. Logic}, 156(2-3):274--286, 2008.

\bibitem[KZ11]{EFAI}
Jonathan Kirby and Boris Zilber.
\newblock Exponential fields and atypical intersections.
\newblock Submitted, arXiv:1108.1075v1, August 2011.

\bibitem[Mac96]{Macintyre96}
Angus~J. Macintyre.
\newblock Exponential algebra.
\newblock In {\em Logic and algebra (Pontignano, 1994)}, volume 180 of {\em
  Lecture Notes in Pure and Appl. Math.}, pages 191--210. Dekker, New York,
  1996.

\bibitem[Zil05]{Zilber05peACF0}
Boris Zilber.
\newblock Pseudo-exponentiation on algebraically closed fields of
  characteristic zero.
\newblock {\em Ann. Pure Appl. Logic}, 132(1):67--95, 2005.

\bibitem[Zil06]{Zilber06covers}
Boris Zilber.
\newblock Covers of the multiplicative group of an algebraically closed field
  of characteristic zero.
\newblock {\em J. London Math. Soc. (2)}, 74(1):41--58, 2006.

\end{thebibliography}

\end{document}